\documentclass[11pt]{article}
\usepackage{authblk}
\usepackage{amsmath}
\usepackage{amssymb}
\usepackage{amsfonts}
\usepackage{comment}
\usepackage{color}
\usepackage{mathtools}
\usepackage{hyperref}
\usepackage{tikz}

\usepackage{amsthm}

\newtheorem*{notation*}{Notation}

\newtheorem{theorem}{Theorem}[section]

\newtheorem{lemma}[theorem]{Lemma}
\newtheorem{claim}[theorem]{Claim}
\newtheorem{observation}[theorem]{Observation}

\newtheorem*{remark*}{Remark}
\numberwithin{equation}{section}

\DeclareMathOperator{\br}{br}
\DeclareMathOperator{\V}{Var}

\theoremstyle{definition}
\newtheorem*{definition}{Definition}

\def \P {\mathbb{P}}
\def \E {\mathbb{E}}

\def \B {\mathcal{B}}
\def \Z {\mathcal{Z}}
\def \M {\mathcal{M}}

\def \Y {\mathcal{Y}}
\def \x {\mathrm{x}}
\def \G {\Gamma}
\def \D {\Delta}

\def \e {\epsilon}

\def \c {\chi_B}

\title{Extremal bipartite independence number and balanced coloring}
\author{Debsoumya Chakraborti\thanks{This work was supported by the Institute for Basic Science (IBS-R029-C1)}}
\affil[]{\small Discrete Mathematics Group (DIMAG), Institute for Basic Science (IBS), Daejeon,~South~Korea}
\affil[]{\small Email: \texttt{debsoumya@ibs.re.kr}}
\begin{document}
\sloppy
\maketitle
\begin{abstract}
In this paper, we establish a couple of results on extremal problems in bipartite graphs. Firstly, we show that every sufficiently large bipartite graph with average degree $D$ and with $n$ vertices on each side has a balanced independent set containing $(1-\e) \frac{\log D}{D} n$ vertices from each side for small $\e > 0$. Secondly, we prove that the vertex set of every sufficiently large balanced bipartite graph with maximum degree at most $\D$ can be partitioned into $(1+\e)\frac{\D}{\log \D}$ balanced independent sets. Both of these results are algorithmic and best possible up to a factor of 2, which might be hard to improve as evidenced by the phenomenon known as `algorithmic barrier' in the literature. The first result improves a recent theorem of Axenovich, Sereni, Snyder, and Weber in a slightly more general setting. The second result improves a theorem of Feige and Kogan about coloring balanced bipartite graphs. 
\end{abstract}

\section{Introduction}

This paper first deals with a bipartite analogue of the Tur\'an's theorem \cite{T41} for complete graphs, which is regarded as a cornerstone of extremal graph theory (see, e.g., \cite{FS} for a survey). Next, we discuss a bipartite analogue of the celebrated Johansson-Molloy Theorem on the chromatic number of a triangle-free graph with a given maximum degree (see, e.g., \cite{B19,M19,MR}). Some seemingly simple problems in the bipartite setting (such as finding the smallest possible `bipartite independence number' of a bipartite graph with maximum degree three) are not yet resolved despite some effort (see, e.g., \cite{ASSW,EMR}). In this paper, we address a few such problems. 

Suppose that we are given a bipartite graph $G=(U \cup V, E)$ with a prescribed vertex bipartition $(U, V)$ and edge set $E$. A balanced bipartite independent set (or bi-hole) of size $t$ in $G$ is a pair $(X,Y)$ where $X \subseteq U$ and $Y \subseteq V$ such that $|X| = |Y| = t$ and there are no edges in $E$ with one endpoint in $X$ and the other in $Y$. The size of the largest bi-hole, referred to as the bipartite independence number, can be viewed as a natural bipartite analogue of the standard independence number. Our first main result is the following. 

\begin{theorem} \label{th1}
For each $\e > 0$, there exists $D_0 = D_0(\e)$ such that the following holds. For each $D \ge D_0$, there is $N_0=N_0(D)$ such that if $G$ is a balanced bipartite graph with average degree $D$ and with $n \ge N_0$ vertices on each side, then $G$ contains a bi-hole of size $(1-\e)\frac{\log D}{D}n$.
\end{theorem}

Theorem~\ref{th1} improves a recent result of Axenovich, Sereni, Snyder, and Weber \cite{ASSW}. They studied the function $f(n,\D)$, which is defined as follows: The function $f(n,\D)$ denotes the largest $k$ such that any bipartite graph $G = (U \cup V, E)$ with $n$ vertices on each of the sides $U$ and $V$, and with maximum degree of $U$ being at most $\D$, contains a bi-hole of size $k$. They determined the correct asymptotic order of $f(n,\D)$ for sufficiently large but fixed $\D$ and growing $n$.

\begin{theorem}[\cite{ASSW}] \label{Axenovich}
For each $0 < \e < 1$, there exists $\D_0 = \D_0(\e)$ such that the following holds. For each $\D \ge \D_0$, there is $N_0 = N_0(\D)$ such that for any $n \ge N_0$, we have that 
\begin{equation*}
    \frac{1}{2} \cdot \frac{\log \D}{\D} \cdot n \le f(n,\D) \le (2+\e) \cdot \frac{\log \D}{\D} \cdot n.
\end{equation*}
\end{theorem}

Note that Theorem~\ref{th1} implies that $f(n,\D)\ge (1-\e) \frac{\log \D}{\D} \cdot n$ because the function $x\rightarrow \frac{\log x}{x}$ is decreasing when $x>0$. We remark that using Theorem \ref{Axenovich}, Ehard, Mohr, and Rautenbach \cite{EMR} gave an easy proof of Theorem~\ref{th1} with a worse bound of $\frac{\log D}{8D}n$. The upper bound in Theorem~\ref{Axenovich} comes from considering the random bipartite graph $G_{n,n,\D/n}$ (the random bipartite graph $G_{n,n,p}$ is a bipartite graph with $n$ vertices on each side where each of the possible $n^2$ edges is present independently with probability $p$). Our proof of Theorem~\ref{th1} yields a randomized algorithm and matches the best bound that can be achieved by an efficient algorithm to find a large bi-hole of $G_{n,n,D/n}$. We elaborate in the concluding remarks why further improving this seems hard. 

We next turn our attention to the bipartite analogue of the standard notion of `chromatic number'. A coloring of the vertices of a balanced bipartite graph $G$ is called `balanced' if each color class induces a bi-hole. The coloring number, $\c(G)$, is the minimum number of colors needed for a balanced coloring of a given balanced bipartite graph $G$ (if such a number does not exist, then define the coloring number to be $\infty$). 

Feige and Kogan \cite{FK} observed that the coloring number of bipartite graphs behaves quite differently from the usual chromatic number of graphs. For example, removing an independent set from a graph never increases its chromatic number. However, removing a bi-hole from a bipartite graph may increase its coloring number. In fact, the remaining graph may not have a balanced coloring at all. This behavior poses some challenges in estimating coloring numbers in general. 

\begin{theorem}[\cite{FK}] \label{Feige}
For each $0 < \e < 1$, there exists $\D_0=\D_0(\e)$ such that the following holds. If $G$ is a balanced bipartite graph with maximum degree $\D \ge \D_0$ and with $n \ge (1+\e)2\D$ vertices on each side, then $$\c(G) \le \frac{20\D}{\e^2\log \D}.$$
\end{theorem}

Our second main result improves the above result and essentially removes the factor of $\frac{20}{\e^2}$. 

\begin{theorem} \label{th2}
For each $\e > 0$, there exists $\D_0=\D_0(\e)$ such that the following holds. For each $\D \ge \D_0$, there is $N_0=N_0(\D)$ such that if $G$ is a balanced bipartite graph with maximum degree $\D$ and with $n \ge N_0$ vertices on each side, then $$\c(G) \le (1+\e)\frac{\D}{\log \D}.$$
\end{theorem}

Our proof of Theorem~\ref{th2} is algorithmic and gives a bound that is best possible up to a factor of 2 (one can easily get a lower bound of $\frac{\D}{(2+\e)\log \D}$ by using Theorem~\ref{Axenovich}). Again, for this coloring problem, our bound matches the best known bound that can be achieved by an efficient algorithm in the random bipartite graph $G_{n,n,\D/n}$.

We observe that one cannot strengthen the bounded maximum degree to a bounded average degree condition in Theorem~\ref{th2}. This can be easily seen from the following fact: If a balanced bipartite graph $G$ with $2n$ vertices contains a vertex $v$ with degree $n$ (i.e., $v$ is connected by an edge with all the vertices from the opposite partition), then $G$ does not have a balanced coloring. 
\medskip

\noindent\textbf{Organization.} This paper is organized in the following way. We start with a list of preliminary tools in the next section that will be helpful throughout the paper. We give a proof of Theorem~\ref{th1} in Section 3 by analyzing a natural randomized algorithm to find a large bi-hole in a given bipartite graph. We next give a more sophisticated randomized algorithm in Section 4 to bound the coloring number of a balanced bipartite graph with bounded maximum degree to prove Theorem~\ref{th2}. This proof uses several technical claims, which will be proved in the subsequent section. Finally, we end with a few concluding remarks in Section 6, where we elaborate on some of the points from the introduction. 

Throughout the paper, we omit the use of floor and ceiling signs for clarity of presentation. For an event $A_n$ that depends on $n$, we say that \textit{$A_n$ occurs `w.h.p.'}, if the probability of $A_n$ tends to one as $n$ tends to infinity.

\section{Preliminaries}

We start with a couple of simple observations regarding balanced coloring from the existing literature, which will be helpful to us later.

\begin{observation} \cite{FK} \label{obs}
A bipartite graph $G=(U\cup V,E)$ has a balanced coloring if and only if the bipartite complement of $G$, i.e., the graph $(U\cup V,\bar{E})$ with $\bar{E}=\{(u,v)\in U\times V:(u,v)\notin E\}$, contains a perfect matching. 
\end{observation}

\begin{lemma} \cite{BHK} \label{easy}
If $G$ is a balanced bipartite graph with maximum degree $\D$ and $n \ge 2\D$ vertices on each side, then $\c(G) \le 2\D + 1$.
\end{lemma}

This lemma gives a weaker upper bound on $\c$ for Theorem~\ref{th2}. Although it appeared in \cite{BHK}, we still prove it to keep our paper self-contained.

\begin{proof}[Proof of Lemma~\ref{easy}]
Let $G$ be a bipartite graph $G$ with maximum degree $\D$ and $n \ge 2\D$ vertices on each side. Consider the bipartite complement $G'$ of $G$. Observe that $G'$ has minimum degree at least $n-\D$. Since $n \ge 2\D$, the minimum degree of $G'$ is at least half of the number of vertices in each side of $G'$. Thus, we deduce that the Hall's conditions hold for $G'$. Thus, by Observation~\ref{obs}, $G$ has a balanced coloring. Now, let $\M = \{e_1, e_2, \ldots, e_n\}$ be a perfect matching of $G'$. We now show that we can greedily color the vertices of $G$ using $2\D + 1$ colors so that both the vertices corresponding to each edge of $\M$ get the same color. Indeed, suppose we have already colored the vertices corresponding to $e_1, e_2, \ldots, e_t$ for some $t < n$. Now, the total number of neighbors of the vertices $u, v$ forming $e_{t+1}$ is at most $2\D$; thus, there must be at least one color left that is used in none of the neighbors of $u$ and $v$. We can use that color for both $u$ and $v$. Thus, each color appears the same number of times in both parts, proving Lemma~\ref{easy}. 
\end{proof}

We need some probabilistic tools to prove Theorems~\ref{th1} and \ref{th2}. We start with a few of the most frequently used probabilistic bounds.

\begin{lemma} [Markov's inequality] \label{markov}
If $X$ is a nonnegative random variable and $t > 0$, then, $$\P[X \ge t] \le \frac{\E(X)}{t}.$$
\end{lemma}

\begin{lemma} [Chebyshev inequality] \label{chebyshev}
If $X$ is a random variable with a finite mean and variance, then, for $t > 0$, $$\P[|X - \E(X)| \ge t] \le \frac{\V(X)}{t^2}.$$
\end{lemma}

We next state the Chernoff bound due to Chernoff \cite{Ch52} and Okamoto \cite{Ok59}. We use the version stated by Janson \cite[Theorem~1]{J02}.

\begin{lemma} [The Chernoff bound] \label{chernoff}
Let $X =\sum_{i=1}^n X_i$, where $X_i$ are independent Bernoulli random variables with $\P[X_i=1]=p_i$. 
Let $\mu=\E(X)=\sum_{i=1}^n p_i$. Then for $t\ge 0$,
\begin{enumerate} 
    \item $\P[X\ge \mu+t]\le e^{-\frac{t^2}{2\mu+2t/3}}$ and
    \item $\P[X\le \mu-t]\le e^{-\frac{t^2}{2\mu}}$.
\end{enumerate}
\end{lemma}

We also need a recent extension \cite{GLSS} of Chernoff bounds to the case when some dependencies between the random variables are allowed. We use the version due to Jukna \cite{J}. To state it, we need the following definition. 

\begin{definition}
A family $Y_1, \ldots, Y_r$ of random variables is \textbf{read-\boldmath$k$} if there exists a sequence $X_1, \ldots, X_m$ of independent random variables, and a sequence $S_1, \ldots, S_r$ of subsets of $[m] = \{1, \ldots, m\}$ such that
\begin{itemize}
    \item each $Y_i$ is some function of $(X_j : j \in S_i)$, and
    \item no element of $[m]$ appears in more than $k$ of the $S_i$'s.
\end{itemize}
\end{definition}

\begin{theorem} [Chernoff bound for dependent random variables, \cite{J}] \label{extension}
Let $Y_1, \ldots, Y_r$ be a family of read-$k$ indicator variables with $\P[Y_i = 1]=p_i$, and let $p$ be the average of $p_1, \ldots, p_r$. Then for any $\e > 0$, $$\P[|(Y_1 + \cdots + Y_r) - pr| \ge \e r] \le 2 e^{-2\e^2r/k}.$$
\end{theorem}

We use the asymmetric version of the local lemma \cite{EL}. We state the version from \cite{AS}.

\begin{lemma} [The local lemma, \cite{AS}] \label{lll}
Let $A_1, \ldots, A_n$ be events in an arbitrary probability space. A directed graph $\mathcal{D}=(V,E)$ on the set of vertices $V=[n]$ is called a dependency digraph for the events $A_1, \ldots, A_n$ if for each $i$, $1 \le i \le n$, the event $A_i$ is mutually independent of all the events $\{A_j : (i,j) \not\in E\}$. Suppose that $\mathcal{D}=(V,E)$ is a dependency digraph for the above events and suppose there are real numbers $x_1, \ldots, x_n$ such that $0 \le x_i < 1$ and $\P[A_i] \le x_i \prod_{(i,j)\in E} (1-x_j)$ for all $1 \le i \le n$. Then, with positive probability no event $A_i$ holds.
\end{lemma}

We want to mention that there are algorithmic versions of the local lemma (see, e.g., \cite{MT,P}). Thus, we can have an efficient randomized algorithm to get the desirable choice of events when we use the local lemma.  We will not further discuss this point inside the proofs.

\section{Finding large bipartite independent sets}

Let $G = (U \cup V, E)$ be an $n$ by $n$ bipartite graph with $|E| = D n$. Let $S_U$ (and $S_V$) denote the set of all vertices in $U$ (and $V$) with degree more than $\frac{D}{\e^2}$. A simple double counting gives us $|S_U|\cdot \frac{D}{\e^2} \le D n$. Thus, we have $|S_U|\le \e^2 n$ and symmetrically, $|S_V|\le \e^2 n$. Consequently, we can remove exactly $\e^2 n$ vertices from both sides to make sure that the maximum degree of the induced graph on the remaining vertices is at most $\frac{D}{\e^2}$. Thus, it is enough to prove Theorem~\ref{th1} with the extra assumption that the maximum degree of the underlying graph is at most $\frac{D}{\e^2}$. This will be crucial in applying certain concentration bounds while analyzing our randomized algorithm. We can assume that $0 < \e < \frac{1}{10}$. Throughout the proof, wherever needed, we will use that $D$ is sufficiently large with respect to $\e$ and $n$ is sufficiently large with respect to $D$. 

The algorithm is straightforward and natural. First, we pick the vertices in $U$ independently with probability $(1-\e/2)\frac{\log D}{D}$. Let $U'$ denote the set of all the vertices picked from $U$. Let $V'$ denote the set of vertices in $V$ that do not have any neighbor in $U'$. To prove Theorem~\ref{th1}, it is enough to show that the sizes of $U'$ and $V'$ are both at least $(1-\e)\frac{\log D}{D}n$ with positive probability. These are shown in the following couple of claims.

\begin{claim} \label{Uside}
W.h.p., we have that $|U'| \ge (1-\e)\frac{\log D}{D}n$.
\end{claim}

\begin{proof}
Let $X_u$ denote the indicator random variable for the event that the vertex $u \in U$ is picked. It is clear that $|U'| = \sum_{u \in U} X_u$. A straightforward application of the Chernoff bound (Lemma~\ref{chernoff}) yields our claim. 
\end{proof}

\begin{claim}
W.h.p., we have that $|V'| \ge (1-\e)\frac{\log D}{D}n$.
\end{claim}

\begin{proof}
For each vertex $v \in V$, let $Y_v$ denote the indicator random variable for the event that no neighbor of $v$ is picked from $U$. It is clear that $|V'| = \sum_{v \in V} Y_v$. We first compute the expected size of $|V'|$. For each $v \in V$, the probability that none of its neighbors are picked is exactly $\left(1-(1-\e/2)\frac{\log D}{D}\right)^{d(v)}$, where $d(v)$ is the degree of $v$. Now, using Jensen's inequality, we have the following.

\begin{align*}
\E(|V'|) = \sum_{v \in V} \left(1-(1-\e/2)\frac{\log D}{D}\right)^{d(v)} 
&\ge n \left(1-(1-\e/2)\frac{\log D}{D}\right)^{D} \\
&\ge n e^{-(1-\e/4)\log D} \\ 
&= \frac{n}{D^{1-\e/4}}.
\end{align*}

We next use Theorem~\ref{extension} to show concentration of the random variable $|V'|$. We claim that the family of random variables $\{Y_v : v \in V\}$ is read-$\frac{D}{\e^2}$. It is clear by observing the following facts.
\begin{itemize}
    \item $X_u, u \in U$ are independent random variables,
    \item for each $v \in V$, $Y_v$ is a function of $(X_u : u \in N(v))$, and
    \item no vertex $u \in U$ is adjacent to more than $\frac{D}{\e^2}$ vertices in $V$.
\end{itemize}

Thus, a straightforward application of Theorem~\ref{extension} on the random variables $Y_v, v \in V$ shows us that $\P[|V'| \le (1-\e)\frac{\log D}{D}n] \le e^{-\Omega_D (n)}$. This finishes the proof of Theorem~\ref{th1}.

\end{proof}

\section{Balanced colorings of bipartite graphs}

In this section, we prove Theorem~\ref{th2} through a series of claims. We later prove these claims in the next section.

\begin{proof}[Proof of Theorem~\ref{th2}]
We can assume that $0 < \e < \frac{1}{10}$. Let $G=(U \cup V, E)$ be an $n$ by $n$ bipartite graph with maximum degree $\D$. Similar to the previous section, wherever needed, we use that $\D$ is sufficiently large with respect to $\e$ and $n$ is sufficiently large with respect to $\D$. Suppose that we are given $(1+\e)\frac{\D}{\log \D}$ colors. To prove Theorem~\ref{th2}, we need to show the existence of a balanced coloring of $G$ using these colors. Fix a set $Q$ of $q = (1+\e/2)\frac{\D}{\log \D}$ colors; there are still $\frac{\e \D}{2\log \D}$ colors outside of $Q$. We first color the vertices in $U$ independently and uniformly at random with the colors in $Q$. We obtain the following fact by a simple application of the Chernoff bound similar to the proof of Claim~\ref{Uside} (we omit the details).

\begin{claim} \label{first}
W.h.p., for every color $c \in Q$, the set of all vertices in $U$ with color $c$, denoted by $U_c$, satisfies that $\frac{n}{q} - \frac{n}{\log n} \le |U_c| \le \frac{n}{q} + \frac{n}{\log n}$. 
\end{claim}

Next, we assign a set $Q_v \subseteq Q$ of available colors to each $v \in V$. Let $C_v$ denote the set of all colors already used by some neighbor of $v$. We set $Q_v = Q \setminus C_v$. For each $v \in V$, we now color it independently and uniformly at random with the colors in $Q_v$ if $Q_v$ is non-empty. If for some $v\in V$, the set $Q_v$ of available colors is empty, then we leave the vertex $v$ uncolored. However, we will show that this does not happen for too many vertices in $V$. Denote by $V_c$ the set of all vertices in $V$ that are colored with $c$.

\begin{claim} \label{second}
\
\begin{enumerate}
\item W.h.p., for every pair of colors $c_1, c_2 \in Q$, we have $\left||V_{c_1}|-|V_{c_2}|\right| \le \frac{n}{\log n}$. 
\item W.h.p., for every color $c \in Q$, we have $|V_c| \ge \left(1-\frac{100}{\e^2 \log^2 \D}\right)\frac{n}{q}$. 
\end{enumerate}
\end{claim}

Our strategy is to finish by coloring all the uncolored vertices in $V$ and recoloring some of the vertices in $U$ and $V$ (to make the coloring balanced) by the remaining $\frac{\e \D}{2\log \D}$ colors which are not in $Q$. To this end, we denote by $S_V$ the set of all the uncolored vertices in $V$. 

Note that Claim~\ref{second}(1) together with the fact that $\sum_{c\in Q} |V_c|\le n$ implies that for every $c\in Q$, we have $|V_c|\le \frac{n}{q}+\frac{n}{\log n}$. By this and Claim~\ref{first}, the following holds for every $c\in Q$. 
\begin{equation} \label{uncolor}
    |U_c| - |V_c| \ge -\frac{2n}{\log n}.
\end{equation}

For every color $c \in Q$, if $|U_c| < |V_c|$, then arbitrarily uncolor some vertices of $V_c$ to make sure that the number of vertices colored with $c$ in both parts is exactly $|U_c|$ (this step is necessary to make sure every color class contains the same number of vertices from $U$ and $V$). Due to \eqref{uncolor}, we have uncolored at most $\frac{2qn}{\log n}$ vertices of $V$, denote by $S_0$ the set of all vertices that got uncolored. Let $Q'$ denote all the colors $c \in Q$ such that $|U_c| \ge |V_c|$. Keep in mind that every color outside of $Q'$ appears the same number of times in $U$ and $V$. 

If the size of $S_V$ is small, then we can greedily finish the coloring as demonstrated next. 
Suppose that $|S_V| \le \frac{n}{\D^2}$. Let $S = S_V \cup S_0$. Then, $|S| \le \frac{n}{\D^2} + \frac{2qn}{\log n} \le \frac{2n}{\D^2}$. We now wish to color all the vertices in $S$ and recolor some vertices of $U$ with a new color $c^*$. More precisely, for every color $c \in Q'$, we will recolor exactly $|U_c| - |V_c|$ vertices of $U_c$ by using $c^*$. Since $|S| \le \frac{2n}{\D^2}$, the number of vertices in $U$ with at least one neighbor in $S$ is at most $\frac{2n}{\D}$, and by Claim~\ref{second}(2), we have $\frac{2n}{\D} < |V_c|$. Then, there are at least $|U_c|-\frac{2n}{\D} > |U_c| - |V_c|$ vertices in $U_c$ that do not have any neighbor in $S$. We next choose exactly $|U_c| - |V_c|$ such vertices of $U_c$ for every $c\in Q'$ and recolor them using $c^*$. Thus, we can successfully color $G$ with $q+1$ colors such that every color class induces a bi-hole.

Thus, from now on, we assume that $|S_V| \ge \frac{n}{\D^2}$. This implies that $\sum_{c\in Q} |V_c|\le n - \frac{n}{\D^2}$. Thus, it follows from Claim~\ref{second}(1) that for every $c\in Q$, we have $|V_c|\le \frac{n}{q}\left(1-\frac{1}{\D^2}\right) +\frac{n}{\log n}$. This fact together with Claims~\ref{first} and~\ref{second}(2) implies the following for every $c\in Q$. 

\begin{equation} \label{local}
    0 \le |U_c|-|V_c| \le \frac{100}{\e^2 \log^2 \D} \cdot \frac{n}{q} + \frac{n}{\log n}.
\end{equation}

Thus, we have $S_0=\emptyset$ and $S_V$ is the set of all uncolored vertices in $V$. 
In this case, we desire to get a set $S_U \subset U$ with the same size as $S_V$ (remember that we want a balanced coloring) such that the maximum degree of the graph induced by $(S_U, S_V)$ is small enough to apply Lemma~\ref{easy} and finish the coloring using the remaining $\frac{\e \D}{2\log \D}$ colors not in $Q$. To achieve this, We start by showing that very few vertices of $U$ have many neighbors in $S_V$.

\begin{claim} \label{lessneighbors}
W.h.p., for every color $c \in Q$, at most $\frac{100n\sqrt{\log \D}}{\e^2\D}$ of the vertices $u$ in $U_c$ satisfy that $u$ has more than $\frac{\D}{\log^{3/2} \D}$ neighbors in $S_V$.
\end{claim}

Suppose now, we fix an instance satisfying all the high probability events. Denote by $U^*_c$ the set of all vertices in $U_c$ with at most $\frac{\D}{\log^{3/2} \D}$ neighbors in $S_V$. By Claims~\ref{first} and \ref{lessneighbors}, we have 

\begin{equation}\label{size of U_c^*}
    |U^*_c| \ge \frac{n}{q} - \frac{n}{\log n} - \frac{100n\sqrt{\log \D}}{\e^2\D} \ge \frac{n\log \D}{2\D}.
\end{equation}

\begin{claim} \label{residue}
There exists $S_U$ consisting of exactly $|U_c|-|V_c|$ many vertices of $U^*_c$ for all $c$, such that the balanced graph induced by $(S_U,S_V)$ has maximum degree at most $\frac{\D}{\log^{3/2} \D}$.
\end{claim}

Finally, it follows from Lemma~\ref{easy} and the fact that $|S_U| = |S_V| \ge \frac{n}{\D^2} > \frac{2\D}{\log^{3/2} \D}$ that there is a balanced coloring of the induced graph on $(S_U,S_V)$ by the remaining $\frac{\e \D}{2\log \D}$ colors that are not used yet. This finishes the proof of Theorem~\ref{th2} modulo the claims.

\end{proof}

\section{Proofs of intermediate claims}

In this section, we complete the proof of Theorem~\ref{th2} by showing the validity of the claims of the last section.

\begin{proof}[Proof of Claim~\ref{second}]
For every color $c \in Q$, let $Z_c$ be the random variable denoting the number of vertices in $V$ with color $c$. Define $Z = \sum_{c \in Q} Z_c$. Observe that $Z = \sum_{v \in V} I_v$, where $I_v$ is the indicator random variable for the set $Q_v$ being non-empty. Hence, 

\begin{equation} \label{nonempty}
    \E(Z) = \sum_{v \in V} \E(I_v) = \sum_{v \in V} \P[Q_v \neq \emptyset].
\end{equation}

For each vertex $v \in V$, the probability that $Q_v$ is empty is the same as the probability that all the colors of $Q$ appear in the neighborhood of $v$. To estimate this probability, consider the following process which essentially describes another way to choose the colors of the vertices in $N(v)\subseteq U$. Start with an empty set $S_0 = \emptyset$, then at each time step $t>0$, we generate a uniformly random color $c_t$ from $Q$ independently of previous choices and define $S_t = S_{t-1} \cup \{c_t\}$ (note that this is a set, hence even if a color comes more than once, it appears only once in $S_t$). Define $T$ to be the random variable that counts the minimum number of time step $t$ such that $|S_t| = q$. Now, observe: 

\begin{equation} \label{relate}
    \P[Q_v = \emptyset] = \P[T \le d(v)] \le \P[T \le \D].
\end{equation}

The random variable $T$ is well-studied and estimating it is known as the `coupon collector's problem' in the literature (see, e.g., \cite{LPW}). To keep our paper self-contained, we estimate the lower tail of $T$ by a simple application of Chebyshev inequality.

\begin{lemma} \label{coupon}
$\P[T \le \D] < \frac{50}{\e^2 \log^2 \D}$.
\end{lemma}

\begin{proof}
For each $1 \le j \le q$, let the random variable $T_j$ denote the first time step $t$ for which $|S_t| = j$ (define $T_0 = 0$). Clearly, $T_q = T$. Note that the random variable $T_j - T_{j-1}$ denotes the time needed for a new color to be added in our collection as $j$-th color. Thus, $T_j - T_{j-1}$ has a geometric distribution with probability $\frac{q-j+1}{q}$. Remember that a random variable with geometric distribution with probability $p$ has expectation $\frac{1}{p}$ and variance $\frac{1-p}{p^2}$. It follows that

\begin{equation} \label{couponexpectation}
    \E(T) = \sum_{j=1}^q \E(T_j - T_{j-1}) = \sum_{j=1}^q \frac{q}{q-j+1} \ge q \int_1^{q+1} \frac{1}{x} dx \ge q \log q.
\end{equation}

Since $q=(1+\e/2)\frac{\D}{\log \D}$, we have $\log q \ge \log \D - \log\log \D \ge \frac{1+\e/4}{1+\e/2}\cdot \log \D$ (where we use the fact that $\D$ is much larger with respect to $\e$). This together with \eqref{couponexpectation} imply the following:

\begin{equation} \label{modcouponexpectation}
    \E(T) \ge (1+\e/4)\D.
\end{equation}

Furthermore, observe that the random variables $T_j - T_{j-1}, j \in [q]$ are independent, and thus, we have the following.

\begin{align} \label{variance}
    \V(T) = \sum_{j=1}^q \V(T_j - T_{j-1}) &\le \sum_{j=1}^q \frac{q^2}{(q-j+1)^2} \nonumber \\
    &\le q^2 \left(1 + \int_1^q \frac{1}{x^2} dx\right) < 2q^2.
\end{align}

Using \eqref{modcouponexpectation}, \eqref{variance}, and Chebyshev inequality (Lemma~\ref{chebyshev}), we have the following.

\begin{align*}
    \P[T \le \D] \le \P\left[T-\E(T) \le -\frac{\e \D}{4}\right] 
    \le \frac{16 \V(T)}{\e^2 \D^2}
    < \frac{50}{\e^2 \log^2 \D}.
\end{align*}
\end{proof}

Thus, using \eqref{nonempty}, \eqref{relate}, and Lemma~\ref{coupon}, we have that $\E(Z) \ge \left(1 - \frac{50}{\e^2 \log^2 \D}\right)n$. By symmetry, $Z_c$ has identical distribution for all $c \in Q$. Thus, by the linearity of expectation, the following holds for every $c \in Q$. \begin{equation}\label{expectation of Z_c}
    \E(Z_c) = \frac{\E(Z)}{q} \ge \left(1 - \frac{50}{\e^2 \log^2 \D}\right)\frac{n}{q}.
\end{equation}

Next, to prove both of the parts of Claim~\ref{second}, we use Theorem~\ref{extension} to show the concentration of each $Z_c$ around its mean.
Fix a color $c \in Q$. For $v \in V$, let $Y_v$ be the indicator random variable for the event that $v$ is colored with $c$. Clearly, $Z_c = \sum_{v \in V} Y_v$. To apply Theorem~\ref{extension}, we wish to show that the family of random variables $\{Y_v : v \in V\}$ is read-$\D$. For $u \in U$, let $X_u$ be the random variable denoting the color chosen for $u$. In order to model the random variables $Y_v$ conveniently, for $v \in V$, let $X'_v$ be independent random variables with continuous uniform distribution on the interval $[0,1)$. For the convenience of our analysis, we now specify how we assign colors to $v \in V$ independently and uniformly at random from the set $Q_v \subseteq Q := [q]$ of available colors. For each $v \in V$, if $Q_v$ is non-empty, then color $v$ with the $j$-th smallest color from $Q_v$, where $j$ satisfies $\frac{j-1}{|Q_v|} \le X'_v < \frac{j}{|Q_v|}$. Now, it is clear that the following facts hold.
\begin{itemize}
    \item $\{X_u : u \in U\} \cup \{X'_v : v \in V\}$ are independent random variables,
    \item for each $v \in V$, the random variable $Y_v$ is a function of $X'_v$ and $(X_u : u \in N(v))$, and
    \item no vertex $u \in U$ is adjacent to more than $\D$ vertices in $V$.
\end{itemize}

Thus, the family of random variables $\{Y_v : v \in V\}$ is read-$\D$. Finally, by applying Theorem~\ref{extension}, the following holds for every $c\in Q$. 
$$\P\left[\left|Z_c - \frac{\E(Z)}{q}\right| \ge \frac{n}{2\log n}\right] = e^{-\Omega_{\D}(n/\log^2 n)}.$$

This together with \eqref{expectation of Z_c} and a simple union bound over all colors finishes the proof of Claim~\ref{second}.
\end{proof}

\begin{proof}[Proof of Claim~\ref{lessneighbors}]
For every color $c \in Q$, let $\Z_c$ be the random variable denoting the number of vertices $u \in U$ with color $c$ and more than $\frac{\D}{\log^{3/2} \D}$ neighbors in $S_V$. Define $\Z = \sum_{c \in Q} \Z_c$. Observe that $\Z = \sum_{u \in U} A_u$, where $A_u$ is the indicator random variable for the event that $u$ has more than $\frac{\D}{\log^{3/2} \D}$ neighbors in $S_V$. For $u \in U$, define the random variable $B_u = \sum_{v \in N(u)} I^c_v$, where $I^c_v$ is the indicator random variable for the set $Q_v$ being empty. Thus, for each $u \in U$, we have that $A_u = 1$ if and only if $B_u > \frac{\D}{\log^{3/2} \D}$. Now, using \eqref{relate} and Lemma~\ref{coupon}, we have the following.

\begin{equation} \label{bu}
    \E(B_u) = \sum_{v \in N(u)} \E(I^c_v) = \sum_{v \in N(u)} \P[Q_v = \emptyset] < \frac{50 \D}{\e^2 \log^2 \D}.
\end{equation}

Thus, by \eqref{bu} and a simple application of Markov's inequality (Lemma~\ref{markov}), we have:

\begin{equation*} 
    \E(A_u) = \P[A_u = 1] = \P\left[B_u > \frac{\D}{\log^{3/2} \D}\right] < \frac{50}{\e^2 \log^{1/2} \D}.
\end{equation*} 

Thus, $\E(\Z) = \sum_{u \in U} \E(A_u) < \frac{50n}{\e^2 \log^{1/2} \D}$. By symmetry, every $\Z_c$ has the same distribution. Hence, by the linearity of expectation, we have that $\E(\Z_c) = \frac{\E(\Z)}{q} < \frac{50n\log^{1/2} \D}{\e^2 \D}$. We next complete the proof of our claim by using Theorem~\ref{extension} to show the concentration of each $\Z_c$ around its mean.

Fix a color $c \in Q$. For $u \in U$, let $\Y_u$ be the indicator random variable for the event that $u$ has color $c$ and $u$ has more than $\frac{\D}{\log^{3/2} \D}$ neighbors in $S_V$. Clearly, $\Z_c = \sum_{u \in U} \Y_u$. We now wish to show that the family of random variables $\{\Y_u : u \in U\}$ is read-$(\D^2 + 1)$. Remember that $X_u$ is the random variable denoting the color of $u \in U$. For convenience, for $u \in U$, define $\G(u)$ to be the set of all vertices in $U$ at distance exactly two from $u$. Now, observe the following:
\begin{itemize}
    \item $\{X_u : u \in U\}$ are independent random variables,
    \item for each $u \in U$, the random variable $\Y_u$ is a function of $X_u$ and $(X_{u'} : u' \in \G(u))$, and
    \item for each $u \in U$, the random variable $X_u$ affects at most $|\G(u)| + 1 \le \D^2 + 1$ many random variables in $\{\Y_u : u \in U\}$.
\end{itemize}

Thus, the family of random variables $\{\Y_u : u \in U\}$ is read-$(\D^2 + 1)$ and a simple application of Theorem~\ref{extension} like before yields Claim~\ref{lessneighbors}.
\end{proof}

\begin{proof}[Proof of Claim~\ref{residue}]
We make use of the local lemma to prove this claim. Include every $u \in U$ independently in a set $S'_U$ with probability $p:=\frac{1}{\log^{7/4} \D}$. For every $v \in S_V$, assign a bad event $B_v$ which denotes that $v$ has more than $\frac{\D}{\log^{3/2} \D}$ neighbors in $S'_U$. For every color $c \in Q$, assign a bad event $A_c$ which denotes that $|S'_U \cap U^*_c| \le \frac{n}{\D \log^{7/8} \D}$. Let us first calculate the probabilities of these bad events. For convenience, denote by $\B(n,p)$ the binomial distribution with the parameters $n$ and $p$. By the Chernoff bound (Lemma~\ref{chernoff}) and \eqref{size of U_c^*}, we obtain the following.

\begin{align}
    \P[B_v] \le \P\left[\B(d(v), p) \ge \frac{\D}{\log^{3/2} \D}\right] \le \P\left[\B(\D, p) \ge \frac{\D}{\log^{3/2} \D}\right] \le e^{-\D^{3/4}}. \label{bv}
\end{align}

\begin{align}
    \P[A_c] \le \P\left[\B\left(|U^*_c|, p\right) \le \frac{n}{\D \log^{7/8} \D}\right] 
    &\le \P\left[\B\left(\frac{n\log \D}{2\D}, p\right) \le \frac{n}{\D \log^{7/8} \D}\right] \nonumber \\
    &\le e^{-\frac{n}{\D \log \D}}. \label{ac}
\end{align}

For $v \in S_V$, let $\G(v)$ denote the set of all vertices in $S_V$ which are in distance exactly 2 from~$v$. Clearly, $|\G(v)| \le \D^2$ for all $v \in S_V$. Note that $B_v$ is mutually independent of all the events $\{B_{v'}:v' \not\in \Gamma(v)\}$. To verify the hypothesis of Lemma~\ref{lll}, set $x_v := e^{-\sqrt{\D}}$ for each $v \in S_V$ and $\x_c := e^{-n/\D^2}$ for each $c \in Q$. We now have the following for each $v \in S_V$.

\begin{align}
    x_v \prod_{v' \in \G(v)} (1 - x_{v'}) \prod_{c \in Q} (1 - \x_c) &\ge e^{-\sqrt{\D}} \left(1 - e^{-\sqrt{\D}}\right)^{\D^2} \left(1 - e^{-n/\D^2}\right)^q \nonumber \\
    &\ge \frac{1}{2} e^{-\sqrt{\D}} \ge \P[B_v], \label{badv}
\end{align}
where in the last step we have used \eqref{bv}. Similarly, we have the following for each $c \in Q$.

\begin{align}
    \x_c \prod_{v \in S_V} (1 - x_v) \prod_{c' \in Q} (1 - \x_{c'}) &\ge e^{-n/\D^2} \left(1 - e^{-\sqrt{\D}}\right)^n \left(1 - e^{-n/\D^2}\right)^q \nonumber \\
    &\ge e^{-n/\D^2} \cdot e^{-n/\D^2} \cdot \frac{1}{2} \ge \P[A_c], \label{badc}
\end{align}
where in the last step we have used \eqref{ac}. Thus, by \eqref{badv}, \eqref{badc}, and using Lemma~\ref{lll}, we have a choice of $S'_U$ such that none of $B_v$ and $A_c$ holds. Now, for each $c \in Q$, choose $|U_c|-|V_c|$ many vertices from $S'_U \cap U^*_c$ and include them in our desirable set $S_U$ (this can be done because of \eqref{local}). It is clear that we still have the property that no vertices in $S_V$ has more than $\frac{\D}{\log^{3/2} \D}$ neighbors in~$S_U$. Remember that for each $c \in Q$, all vertices in $U^*_c$ have at most $\frac{\D}{\log^{3/2} \D}$ neighbors in $S_V$. This proves Claim~\ref{residue}.
\end{proof}

This finishes the proof of Theorem~\ref{th2}.

\section{Concluding remarks}

How good is the estimate of Lemma~\ref{coupon}? There are some classical `central limit theorem' type results on coupon collector's problem (see, e.g., \cite{ER,LPW}), which do not seem to help us in improving Lemma~\ref{coupon}. However, if one uses a recent result (Theorem~1.9.3 in \cite{D}), then it seems possible to prove Theorem~\ref{th2} avoiding Claims~\ref{lessneighbors} and \ref{residue} (thus, we would not need the local lemma). Nevertheless, we refrain from using such a strong result and keep our paper self-contained. 

We remark that finding the largest bi-hole of a bipartite graph is an NP-hard problem. To see this and some inapproximability results on the bipartite independence number, the interested readers can have a look at \cite{FK04}. Naturally, one can expect the problem of finding the coloring number of a bipartite graph to be even more challenging.

We next discuss why the current known upper bound of Theorem~\ref{th1} and lower bound of Theorem~\ref{th2} can be hard to improve by considering the appropriate random bipartite graphs. To show the upper bound of Theorem~\ref{Axenovich}, the authors \cite{ASSW} essentially proved that the random bipartite graph $G_{n,n,D/n}$ cannot have a bi-hole of size $(2+\e) \frac{\log D}{D} n$ w.h.p. It can be shown (using essentially the same arguments as in \cite{F90} or \cite{FK15}) that this upper bound is asymptotically tight for the bipartite independence number of $G_{n,n,D/n}$ w.h.p. Thus, by considering random bipartite graphs, it is not possible to obtain a bi-hole of size $(2+\e) \frac{\log D}{D} n$ in Theorem~\ref{th1} . It can also be shown by a standard argument (similar to the one for the chromatic number of the random graph $G_{n,\D/n}$; see, e.g., \cite{FK15}) that the coloring number of the random bipartite graph $G_{n,n,\D/n}$ is concentrated around $\frac{\D}{2\log \D}$ w.h.p. Thus, by considering random bipartite graphs, the upper bound on $\c(G)$ in Theorem~\ref{th2} cannot be improved by a factor more than 2.

We next reason why we believe that improving the current gap of a factor of $2$ between lower and upper bounds in Theorems~\ref{th1} and \ref{th2} can be challenging. Before discussing it, we mention the situation for a similar problem in graphs (not restricted to bipartite graphs). The best known lower and upper bounds for the largest possible chromatic number of a triangle-free graph with a bounded maximum degree have a multiplicative gap of two. However, it is believed to be hard to improve this gap (see, e.g., \cite{AC,M19,ZK}). We experience a similar situation in the bipartite setting, as demonstrated next.

A simple greedy algorithm obtains a bi-hole of size $(1-\e)\frac{\log D}{D}n$ in the random bipartite graph $G_{n,n,D/n}$ w.h.p. (e.g., the same method as in Exercise 6.7.20 of \cite{FK15} works here). However, no efficient (polynomial time) algorithm (deterministic or randomized) is known to find a significantly larger bi-hole (see, e.g., \cite{AC,ZK}). This shows some difficulty of improving Theorem~\ref{th1}, it seems especially challenging to find an efficient algorithm to find a significantly larger bi-hole in Theorem~\ref{th1} (because, an algorithm for Theorem~\ref{th1}  will likely find a similar-sized bi-hole in $G_{n,n,D/n}$). On the other hand, since there is no efficient algorithm known to find a bi-hole in $G_{n,n,\D/n}$ of size significantly larger than $\frac{\log \D}{\D}n$, we do not have any efficient algorithm to color $G_{n,n,\D/n}$ using significantly less than $\frac{\D}{\log \D}$ colors. Our bound of Theorem~\ref{th2} matches this and extends this to efficiently color any bipartite graph with maximum degree $\D$ with about $\frac{\D}{\log \D}$ colors.

We next briefly discuss some related problems to Theorem~\ref{th1} in the literature. We would suggest the readers have a look at Section 2 of \cite{ASSW} to see a more detailed description of various connections with Theorem~\ref{th1} or \ref{Axenovich}. As mentioned in \cite{ASSW}, they are related to the bipartite version of the Erd\H{o}s-Hajnal conjecture (see, e.g., \cite{ATW,EHP}), the bipartite Ramsey numbers (see, e.g., \cite{CR,C08}), and the Zarankiewicz function (see, e.g., \cite{BGMV07,BGMV08,FS,GO,GST}). To see the connection with the bipartite Ramsey number, for bipartite graphs $H_1$ and $H_2$, let the bipartite Ramsey number $\br(H_1, H_2)$ be the smallest $N$ such that any red-blue edge-coloring of the complete bipartite graph $K_{N,N}$ contains either a red copy of $H_1$ or a blue copy of $H_2$. For results on this topic, see, e.g., Beineke and Schwenk \cite{BS}, Caro and Rousseau \cite{CR}, Conlon \cite{C08}, Hattingh and Henning \cite{HH}, Irving \cite{I78}, Lin and Li \cite{LL}, and Thomason \cite{T82}. As an application of Theorem~\ref{th1}, we obtain that $\br(K_{1,\D},K_{n,n}) \lesssim \frac{\D}{\log \D} n$ for sufficiently large but fixed $\D$ and growing $n$. 

We end by suggesting two directions for future research. Firstly, it will be interesting to study multi-partite analogues of Theorems~\ref{th1} and~\ref{th2}. For example, one can define `tri-hole' in a tripartite graph as an independent set with the same number of vertices in all three parts. It might be worth estimating the size of the largest tri-hole in a tripartite graph with a bounded average degree or a bounded local degree. The straightforward extensions of the methods used in this paper do not seem to work for $k$-partite graphs when $k \ge 3$.

There is a recent result by Kogan \cite{K} on a generalization of the notion of bipartite independence number. They bounded the largest $k$ for which a given $n$ by $n$ bipartite graph has a $k$ by $k$ induced $d$-degenerate subgraph. This can be studied in the context of Theorem~\ref{th1}. For example, it is worth investigating if one can improve the trivial bound obtained by Theorem~\ref{th1} to get a significantly larger balanced $d$-degenerate subgraph.  

\section*{Acknowledgements}

We are thankful to Rutger Campbell, Sang-il Oum, and the anonymous referees for helping us to improve the writing of this paper.

\end{document}